\documentclass[12pt]{amsart}
\usepackage{amsmath} 
\usepackage{amsfonts} 
\usepackage{amssymb} 
\usepackage{amsthm} 
\usepackage{amscd} 
\usepackage{enumerate}
\usepackage{graphicx,xcolor}
\numberwithin{equation}{section}

\theoremstyle{plain}
\newtheorem{theorem}{Theorem}[section]
\newtheorem*{theorem*}{Theorem} 

\newtheorem*{lemma*}{Lemma} 

\newtheorem*{corollary*}{Corollary} 

\newtheorem*{consequence*}{Consequence} 
\newtheorem{proposition}{Proposition}[section]
\newtheorem*{proposition*}{Proposition} 

\newtheorem*{conjecture*}{Conjecture} 

\theoremstyle{definition}
\newtheorem{definition}{Definition}[section]
\newtheorem*{definition*}{Definition} 

\newtheorem*{remark*}{Remark} 
\newtheorem*{remarks*}{Remarks} 

\newtheorem*{question*}{Question} 
\newtheorem*{questions*}{Questions} 

\newtheorem*{example*}{Example} 
\newtheorem*{examples*}{Examples} 

\newtheorem*{exercise*}{Exercise} 
\newtheorem*{exercises*}{Exercises} 

\theoremstyle{plain}

\newtheorem*{thm*}{Théorème} 

\newtheorem*{lemme*}{Lemme} 

\newtheorem*{cor*}{Corollaire} 

\newtheorem*{csq*}{Conséquence} 

\newtheorem*{prop*}{Proposition} 

\newtheorem*{conj*}{Conjecture} 

\theoremstyle{definition}

\newtheorem*{déf*}{Définition} 

\newtheorem*{rem*}{Remarque} 
\newtheorem*{rems*}{Remarques} 

\newtheorem*{ex*}{Exemple} 
\newtheorem*{exs*}{Exemples} 

\newtheorem*{exo*}{Exercice} 
\newtheorem*{exos*}{Exercices} 

\theoremstyle{remark}

\newcommand{\disp}{\displaystyle}

\renewcommand{\.}{{}_{\!}} 

\renewcommand{\a}{\alpha}

\renewcommand{\c}{\gamma}

\renewcommand{\t}{\theta}

\newcommand{\cC}{\mathcal{C}}

\newcommand{\cE}{\mathcal{E}}

\newcommand{\as}{\ \mbox{\raisebox{.085ex}{$:$}\!$=$} \ } 

\newcommand{\st}{~|~} 

\newcommand{\llist}[3]{#1_{#2} , \ldots , #1_{#3}} 
\newcommand{\lset}[3]{\left\{ #1_{#2} , \ldots , #1_{#3} \right\}} 

\newcommand{\inc}{\subseteq} 
\newcommand{\setmin}{\raisebox{0.45ex}{\scriptsize $\smallsetminus$}} 

\renewcommand{\to}{\longrightarrow} 
 
\newcommand{\seqN}[2]{\left( #1_{#2} \right)_{\! #2 \in \NN}} 

\newcommand{\card}[1]{\mathrm{card} \! \left( #1 \right)} 

\newcommand{\NN}{\mathbf{N}} 
\newcommand{\RR}{\mathbf{R}} 

\newcommand{\Rn}[1]{\mathbf{R}^{\! #1 \!}} 

\newcommand{\intr}[1]{\overset{\; _{\circ}}{#1}} 
\newcommand{\clos}[1]{\overline{#1}} 

\renewcommand{\leq}{\leqslant} 
\renewcommand{\geq}{\geqslant} 

\let\oldint\int
\renewcommand{\int}[4]{\oldint_{\! #1}^{#2} \!\!\!\! #3 \mathrm{d} #4} 

\newcommand{\goes}{\rightarrow} 

\let\olddet\det
\renewcommand{\det}[1]{\olddet{\! \left( #1 \right)}} 
 
\newcommand{\norm}[1]{\left\| #1 \right\|} 
 
\newcommand{\bC}{\partial \cC} 

\newcommand{\dC}{d_{\cC}} 

\renewcommand{\inc}{\subseteq} 

\begin{document}

\title[]{Hilbert domains quasi-isometric to normed vector spaces}

\author{Bruno Colbois}
\address{Bruno Colbois, 
Universit\'{e} de Neuch\^{a}tel, 
Institut de math\'{e}matique, 
Rue \'{E}mile Argand 11, 
Case postale 158, 
CH--2009 Neuch\^{a}tel, 
Switzerland}
\email{bruno.colbois@unine.ch}

\author{Patrick Verovic}
\address{Patrick Verovic, 
UMR 5127 du CNRS \& Universit\'{e} de Savoie, 
Laboratoire de math\'{e}matique, 
Campus scientifique, 
73376 Le Bourget-du-Lac Cedex, 
France}
\email{verovic@univ-savoie.fr}

\date{\today}
\subjclass[2000]{Primary: global Finsler geometry, Secondary: convexity}


\begin{abstract}

We prove that a Hilbert domain which is quasi-isometric to a normed vector space is actually a convex polytope. 

\end{abstract}

\maketitle

\bigskip
\bigskip


\section{Introduction} 

A \emph{Hilbert domain} in $\Rn{m}$ is a metric space $(\cC , \dC)$, where $\cC$ is an 
\emph{open bounded convex} set in $\Rn{m}$ and $\dC$ is the distance function on $\cC$ 
--- called the \emph{Hilbert metric} --- defined as follows. 
   
\medskip
   
Given two distinct points $p$ and $q$ in $\cC$, 
let $a$ and $b$ be the intersection points of the straight line defined by $p$ and $q$ 
with $\bC$ so that $p = (1 - s) a + s b$ and $q = (1 - t) a + t b$ with $0 < s < t < 1$. 
Then 
$$
\dC(p , q) \as \frac{1}{2} \ln{\! [a , p , q , b]},
$$ 
where 
$$
[a , p , q , b] \as \frac{1 - s}{s} \times \frac{t}{1 - t} > 1
$$ 
is the cross ratio of the $4$-tuple of ordered collinear points $(a , p , q , b)$. \\ 
We complete the definition by setting $\dC(p , p) \as 0$. 

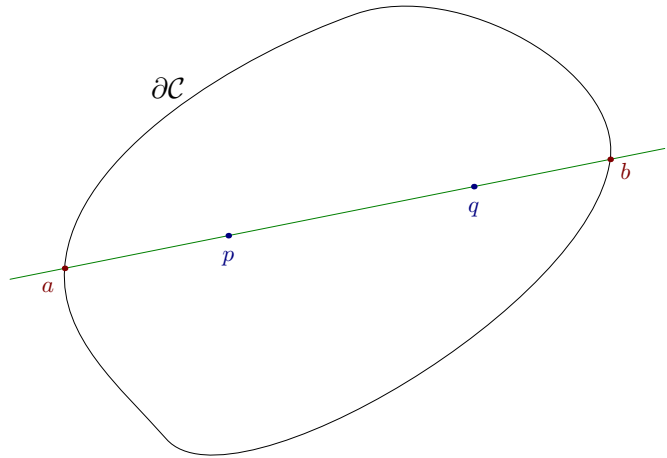
\begin{figure}[h]
   \includegraphics[width=9.7cm,height=6cm,keepaspectratio=true]{fig-1.eps} 
\end{figure}

\bigskip

The metric space $(\cC , \dC)$ thus obtained is a complete non-compact geodesic metric space 
whose topology is the one induced by the canonical topology of $\Rn{m}$ and in which the affine open segments 
joining two points of the boundary $\bC$ are geodesics that are isometric to $(\RR , |\cdot|)$. 

\bigskip

For further information about Hilbert geometry, we refer to \cite{Bus55,BusKel53,Egl97,Gol88} 
and the excellent introduction \cite{Soc00} by Soci\'{e}-M\'{e}thou. 

\bigskip

The two fundamental examples of Hilbert domains $(\cC , \dC)$ in $\Rn{m}$ correspond to the case 
when $\, \cC$ is an ellipsoid, which gives the Klein model of $m$-dimensional hyperbolic geometry 
(see for example \cite[first chapter]{Soc00}), and the case when the closure $\, \clos{\cC} \,$ 
is a $m$-simplex for which there exists a norm $\norm{\cdot}_{\cC}$ on $\Rn{m}$ such that $(\cC , \dC)$ 
is isometric to the normed vector space $(\Rn{m} , \norm{\cdot}_{\cC})$ 
(see \cite[pages 110--113]{dLH93} or \cite[pages 22--23]{Nus88}). 

\bigskip

Much has been done to study the similarities between Hilbert and hyperbolic geometries 
(see for example \cite{ColVero04}, \cite{Ver05} or \cite{Ben06}), 
but little literature deals with the question of knowing to what extend 
a Hilbert geometry is close to that of a normed vector space. 
So let us mention three results in this latter direction which are relevant for our present work. 

\begin{theorem}[\cite{FörKar05}, Theorem~2] \label{thm:isometric} 
   A Hilbert domain $(\cC , \dC)$ in $\Rn{m}$ is isometric to a normed vector space 
   if and only if $\; \cC$ is the interior of a $m$-simplex. 
\end{theorem}

\smallskip

\begin{theorem}[\cite{CVV08}, Theorem~3.1] \label{thm:Lipschitz-equ-2} 
   If $\: \cC$ is an open convex polygonal set in $\Rn{2}$, 
   then $(\cC , \dC)$ is Lipschitz equivalent to Euclidean plane. 
\end{theorem}

\smallskip

\begin{theorem}[\cite{Ber08}, Theorem~1.1. See also \cite{Ver08}] \label{thm:Lipschitz-equ-m} 
   If $\: \cC$ is an open set in $\Rn{m}$ whose closure $\, \clos{\cC}$ is a convex polytope, 
   then $(\cC , \dC)$ is Lipschitz equivalent to Euclidean $m$-space. 
\end{theorem}

\medskip

Recall that a convex \emph{polytope} in $\Rn{m}$ (called a convex \emph{polygon} when $m \as 2$) 
is the convex hull of a finite set of points whose affine span is the whole space $\Rn{m}$. 

\bigskip

In light of these three results, it is natural to ask whether the converse 
of Theorem~\ref{thm:Lipschitz-equ-m} 
--- which generalizes Theorem~\ref{thm:Lipschitz-equ-2} in higher dimensions --- holds. 
In other words, if a Hilbert domain $(\cC , \dC)$ in $\Rn{m}$ 
is quasi-isometric to a normed vector space, what can be said about $\cC$? 
Here, by \emph{quasi-isometric} we mean the following (see \cite{BBI01}): 

\begin{definition} 
   Given real numbers $A \geq 1$ and $B \geq 0$, a metric space $(S , d)$ 
   is said to be $(A , B)$-quasi-isometric to a normed vector space $(V , \norm{\cdot})$ 
   if and only if there exists a map $f : S \to V$ such that 
   $$
   \frac{1}{A} d(p , q) - B \leq \norm{f(p) - f(q)} \leq A d(p , q) + B
   $$ 
   for all $p , q \in S$. 
\end{definition}

\bigskip

We can now state the result of this paper which asserts that the converse 
of Theorem~\ref{thm:Lipschitz-equ-m} is actually true: 

\begin{theorem} \label{thm:main} 
   If a Hilbert domain $(\cC , \dC)$ in $\Rn{m}$ is $(A , B)$-quasi-isometric 
   to a normed vector space $(V , \norm{\cdot})$ for some real constants $A \geq 1$ and $B \geq 0$, 
   then $\, \cC$ is the interior of a convex polytope. 
\end{theorem}

\bigskip
\bigskip
\bigskip


\section{Proof of Theorem \ref{thm:main}} \label{sec:proof} 

The proof of Theorem~\ref{thm:main} is based on an idea developed by F{\"o}rtsch and Karlsson 
in their paper~\cite{FörKar05}. 

\medskip

It needs the following fact due to Karlsson and Noskov: 

\begin{theorem}[\cite{KarNos02}, Theorem~5.2] \label{thm:Gromov-product} 
   Let $(\cC , \dC)$ be a Hilbert domain in $\Rn{m}$ and $x , y \in \bC$ such that $[x , y] \not \inc \bC$. 
   Then, given a point $p_{0} \in \cC$, there exists a constant $K(p_{0} , x , y) > 0$ 
   such that for any sequences $\seqN{x}{n}$ and $\seqN{y}{n}$ in $\, \cC$ 
   that converge respectively to $x$ and $y$ in $\Rn{m}$ 
   one can find an integer $n_{0} \in \NN$ for which we have 
   $$
   \dC(x_{n} , y_{n}) \geq \dC(x_{n} , p_{0}) + \dC(y_{n} , p_{0}) - K(p_{0} , x , y)
   $$ 
   for all $n \geq n_{0}$. 
\end{theorem}

\bigskip

Now, here is the key result which gives the proof of Theorem~\ref{thm:main}: 

\begin{proposition} \label{prop:segments-boundary} 
   Let $(\cC , \dC)$ be a Hilbert domain in $\Rn{m}$ which is $(A , B)$-quasi-isometric 
   to a normed vector space $(V , \norm{\cdot})$ for some real constants $A \geq 1$ and $B \geq 0$. 
   
   Then, if $N = N(A , \norm{\cdot})$ denotes the maximum number of points in the ball 
   $\{ v \in V \st \norm{v} \leq 2 A \}$ whose pairwise distances with respect to $\norm{\cdot}$ 
   are greater than or equal to $1 / (2 A)$, and if $X \inc \bC$ is such that 
   $[x , y] \not \inc \bC$ for all $x , y  \in X$ with $x \neq y$, we have 
   $$
   \card{X} \leq N.
   $$ 
\end{proposition}

\medskip

\begin{proof}~\\ 
Let $f : \cC \to V$ such that 
\begin{equation} \label{equ:quasi-isometry} 
   \frac{1}{A} \dC(p , q) - B \leq \norm{f(p) - f(q)} \leq A \dC(p , q) + B 
\end{equation} 
for all $p , q \in \cC$. 

\medskip

First of all, up to translations, we may assume that $0 \in \cC$ and $f(0) = 0$. 

\medskip

Then suppose that there exists a subset $X$ of the boundary $\bC$ such that 
$[x , y] \not \inc \bC$ for all $x , y  \in X$ with $x \neq y$ and 
$\card{X} \geq N + 1$. So, pick $N + 1$ distinct points $\llist{x}{1}{N + 1}$ in $X$, 
and for each $k \in \{ 1 , \ldots , N + 1 \}$, let $\c_{k} : [0 , +\infty) \to \cC$ 
be a geodesic of $(\cC , \dC)$ that satisfies 
$\c_{k}(0) = 0$, $\disp \lim_{t \goes +\infty} \c_{k}(t) = x_{k}$ in $\Rn{m}$ 
and $\dC(0 , \c_{k}(t)) = t$ for all $t \geq 0$. 

This implies that for all integers $n \geq 1$ and every $k \in \{ 1 , \ldots , N + 1 \}$, we have 
\begin{equation} \label{equ:less} 
   \norm{\frac{f(\c_{k}(n))}{n}} \leq A + \frac{B}{n} 
\end{equation} 
from the second inequality in Equation~\ref{equ:quasi-isometry} with $p \as \c_{k}(n)$ and $q \as 0$. 

\smallskip

On the other hand, Theorem~\ref{thm:Gromov-product} yields the existence of some integer $n_{0} \geq 1$ 
such that 
$$
\dC(\c_{i}(n) , \c_{j}(n)) \geq 2 n - K(0 , x_{i} , x_{j})
$$ 
for all integers $n \geq n_{0}$ and every $i , j \in \{ 1 , \ldots , N + 1 \}$ with $i \neq j$, 
and hence 
\begin{equation} \label{equ:greater} 
   \norm{\frac{f(\c_{i}(n))}{n} - \frac{f(\c_{j}(n))}{n}} 
   \geq 
   \frac{2}{A} - \frac{1}{n} \!\! \left( \! \frac{K(0 , x_{i} , x_{j})}{A} + B \! \right) 
\end{equation} 
from the first inequality in Equation~\ref{equ:quasi-isometry} with $p \as \c_{i}(n)$ and $q \as \c_{j}(n)$. 

\smallskip

Now, fixing an integer 
$n \geq n_{0} + A B + \max{\! \{ K(0 , x_{i} , x_{j}) \st i , j \in \{ 1 , \ldots , N + 1 \} \}}$, 
we get 
$$
\norm{\frac{f(\c_{k}(n))}{n}} \leq 2 A
$$ 
for all $k \in \{ 1 , \ldots , N + 1 \}$ by Equation~\ref{equ:less} together with 
$$
\norm{\frac{f(\c_{i}(n))}{n} - \frac{f(\c_{j}(n))}{n}} \geq \frac{1}{2 A}
$$ 
for all $i , j \in \{ 1 , \ldots , N + 1 \}$ with $i \neq j$ by Equation~\ref{equ:greater}. 

\smallskip

But this contradicts the definition of $N = N(A , \norm{\cdot})$. 

\smallskip

Therefore, Proposition~\ref{prop:segments-boundary} is proved. 
\end{proof}

\bigskip

\begin{remark*} 
Given $v \in V$ such that $\norm{v} = 2 A$, we have $\norm{-v} = 2 A$ 
and $\norm{v - (-v)} = 2 \norm{v} = 4 A \geq 1 / (2 A)$, which shows that $N \geq 2$. 
\end{remark*}

\bigskip

The second ingredient we will need for the proof of Theorem~\ref{thm:main} is the following: 

\begin{proposition} \label{prop:convex-polygon} 
   Let $\, \cC$ be an open bounded convex set in $\Rn{2}$. 
   
   If there exists a non-empty finite subset $Y \!\!$ of the boundary $\, \bC$ such that for every $x \in \bC$ 
   one can find $y \in Y \!\!$ with $[x , y] \inc \bC$, then the closure $\, \clos{\cC}$ is a convex polygon. 
\end{proposition}

\medskip

\begin{proof}~\\ 
Assume $0 \in \cC$ and let us consider the continuous map $\pi : \RR \to \bC$ 
which assigns to each $\t \in \RR$ the unique intersection point $\pi(\t)$ of $\, \bC$ 
with the half-line $\RR_{+}^{*} (\cos{\t} , \sin{\t})$. 

\smallskip

For each pair $(x_{1} , x_{2}) \in \bC \times \bC$, denote by $A(x_{1} , x_{2}) \inc \bC$ the arc segment 
defined by $A(x_{1} , x_{2}) \as \pi([\t_{1} , \t_{2}])$, 
where $\t_{1}$ and $\t_{2}$ are the unique real numbers such that 
$\pi(\t_{1}) = x_{1}$ and $\pi(\t_{2}) = x_{2}$ 
with $\t_{1} \in [0 , 2 \pi)$ and $\t_{1} \leq \t_{2} < \t_{1} + 2 \pi$. 

\medskip

Before proving Proposition~\ref{prop:convex-polygon}, notice that adding a point of $\bC$ to $Y \!$ 
does not change $Y$'s property at all, and therefore we may assume that $\card{Y} \geq 2$. 

\smallskip

So, write $Y \! = \lset{x}{1}{n}$ with $x_{1} = \pi(\t_{1}), \ldots , x_{n} = \pi(\t_{n})$, 
where $\t_{1} \in [0 , 2 \pi)$ and $\t_{1} < \cdots < \t_{n} < \t_{n + 1} \as \t_{1} + 2 \pi$, 
and let $x_{n + 1} \as \pi(\t_{n + 1}) = x_{1}$. 

\medskip

Fix $k \in \{ 1 , \ldots , n \}$ and pick an arbitrary 
$x \in A(x_{k} , x_{k + 1}) \setmin \{ x_{k} , x_{k + 1} \}$. 

\medskip

By hypothesis, one can find $y \in Y \!$ with $[x , y] \inc \bC$.  

\smallskip

Then the convex set $\cC$ is contained in one of the two open half-planes in $\Rn{2}$ 
bounded by the line passing through the points $x$ and $y$, 
and hence either $A(x , y) = [x , y]$, or $A(y , x) = [x , y]$. 

\smallskip

Since $x_{k} \in A(y , x)$ and $x_{k + 1} \in A(x , y)$, we then have 
$x_{k} \in [x , y]$ or $x_{k + 1} \in [x , y]$, 
which yields $A(x_{k} , x) = [x_{k} , x]$ or $A(x , x_{k + 1}) = [x , x_{k + 1}]$. 

\smallskip

Conclusion: $A(x_{k} , x_{k + 1}) = S_{k} \cup S_{k + 1}$, 
where $S_{k} \as \{ x \in A(x_{k} , x_{k + 1}) \st A(x_{k} , x) = [x_{k} , x] \}$ 
and $S_{k + 1} \as \{ x \in A(x_{k} , x_{k + 1}) \st A(x , x_{k + 1}) = [x , x_{k + 1}] \}$. 

\smallskip

Now, the set $S_{k}$ (resp.~$S_{k + 1}$) satisfies $[x_{k} , x] \inc S_{k}$ 
(resp.~$[x , x_{k + 1}] \inc S_{k + 1}$) whenever $x \in S_{k}$ (resp.~$x \in S_{k + 1}$). 

\smallskip

So, if we consider 
$\a_{0} \as \max{\! \{ \t \in [\t_{k} , \t_{k + 1}] \st A(x_{k} , \pi(\t)) = [x_{k} , \pi(\t)] \}}$, 
we have $S_{k} = [x_{k} , \pi(\a_{0})]$ and $S_{k + 1} = [\pi(\a_{0}) , x_{k + 1}]$. 

\smallskip

Hence, $A(x_{k} , x_{k + 1})$ is the union of the two affine segments 
$[x_{k} , \pi(\a_{0})]$ and $[\pi(\a_{0}) , x_{k + 1}]$. 

\medskip

Finally, since $\disp \bC = \bigcup_{k = 1}^{n} A(x_{k} , x_{k + 1})$, 
this implies that $\bC$ is the union of $2 n$ affine segments in $\Rn{2}$, 
and thus $\, \clos{\cC} \,$ is a convex polygon. 
\end{proof}

\bigskip

Before proving Theorem~\ref{thm:main}, let us recall the following useful result, 
where a convex \emph{polyhedron} in $\Rn{m}$ is the intersection of a finite number of closed half-spaces: 

\begin{theorem}[\cite{Kle59}, Theorem~4.7] \label{thm:plane-sections} 
   Let $\, P$ be a convex set in $\Rn{m}$ and $p \in \. \intr{P}$. 
   
   Then $P$ is a convex polyhedron if and only if all its plane sections containing $p$ are convex polyhedra. 
\end{theorem}

\bigskip

\begin{proof}[Proof of Theorem~\ref{thm:main}]~\\ 
Let $(\cC , \dC)$ be a non-empty Hilbert domain in $\Rn{m}$ that is $(A , B)$-quasi-isometric 
to a normed vector space $(V , \norm{\cdot})$ for some real constants $A \geq 1$ and $B \geq 0$. 

\smallskip

According to Theorem~\ref{thm:plane-sections}, it suffices to prove Theorem~\ref{thm:main} for $m \as 2$ 
since any plane section of $\,\. \cC$ gives rise to a $2$-dimensional Hilbert domain 
which is also $(A , B)$-quasi-isometric to $(V , \norm{\cdot})$. 

\medskip

So, let $m \as 2$, and consider the set 
$\cE \as \{ X \inc \bC \st [x , y] \not \inc \bC \ \mbox{for all} \ x , y \in X \ \mbox{with} \ x \neq y \}$. 

\smallskip

It is not empty since $\{ x , y \} \in \cE$ for some $x , y \in \bC$ with $x \neq y$ 
(indeed, $\cC$ is a non-empty open set in $\Rn{2}$), 
which implies together with Proposition~\ref{prop:segments-boundary} that 
$n \as \max{\! \{ \card{X} \st X \in \cE \}}$ does exist and satisfies $2 \leq n \leq N$ 
(recall that $N \geq 2$). 

\medskip

Then pick $Y \! \in \cE$ such that $\card{Y} = n$, write $Y \! = \lset{x}{1}{n}$, 
and prove that for every $x \in \bC$ one can find $k \in \{ 1 , \ldots , n \}$ 
such that $[x , x_{k}] \inc \bC$. 

\smallskip

Owing to Proposition~\ref{prop:convex-polygon}, this will show that $\clos{\cC}$ is a convex polygon. 

\medskip

So, suppose that there exists $x_{0} \in \bC$ satisfying $[x_{0} , x_{k}] \not \inc \bC$ 
for all $k \in \{ 1 , \ldots , n \}$, 
and let us find a contradiction by considering $Z \as Y \cup \{ x_{0} \}$. 

\smallskip

First, since $x_{0} \neq x_{k}$ for all $k \in \{ 1 , \ldots , n \}$ 
(if not, we would get an index $k \in \{ 1 , \ldots , n \}$ such that $[x_{0} , x_{k}] = \{ x_{0} \} \inc \bC$, 
which is false), we have $x_{0} \not \in Y \!$. Hence $\card{Z} = n + 1$. 

\smallskip

Next, since $Y \! \in \cE$ and $[x_{0} , x_{k}] \not \inc \bC$ for all $k \in \{ 1 , \ldots , n \}$, 
we have $Z \in \cE$. 

\smallskip

Therefore, the assumption of the existence of $x_{0}$ yields a set $Z \in \cE$ whose cardinality 
is greater than that of $Y \!$, which contradicts the very definition of $Y \!$. 

\medskip

Conclusion: $\clos{\cC}$ is a convex polygon, and this proves Theorem~\ref{thm:main}. 
\end{proof}

\bigskip
\bigskip
\bigskip


\bibliographystyle{acm}
\bibliography{math-biblio}

\end{document}